\documentclass{amsart}
\usepackage{amsmath, amsthm}

\newcommand\CC{\mathbb C} 
\newcommand\RR{\mathbb R}
\newcommand\Cn{\CC^n}

\newcommand\HH{\mathbb H}
\newcommand\db{\overline{\partial}}  
\newcommand\z{\zeta}
\newcommand\OO{\Omega}
\newcommand\bO{\partial\Omega}
\newcommand\cO{\overline\Omega}
\newcommand{\dzk}[1]{\frac{\partial#1}{\partial\bar{\zeta}_k}}
\newcommand{\dd}[2]{\frac{\partial#1}{\partial#2}} 
\newcommand\BM{Bochner-Martinelli}
\newcommand\pr{$\psi$-regular}

\theoremstyle{theorem}
\newtheorem{theorem}{Theorem}
\newtheorem{corollary}[theorem]{Corollary}
\newtheorem{proposition}[theorem]{Proposition}

\theoremstyle{remark}
\newtheorem{remark}{Remark}

\begin{document}

\title{Quaternionic regularity\break and the $\db$-Neumann problem in $\CC^2$}

\author[A. Perotti]{Alessandro Perotti}
\address{Department of Mathematics,  University of Trento, Via Sommarive 14, 38050 Povo-Trento, ITALY} \email{perotti@science.unitn.it}
\urladdr{http://www.science.unitn.it/\textasciitilde perotti/}

\begin{abstract}
Let $\OO$ be a domain in the quaternionic space $\HH$. We prove a differential criterion that characterizes Fueter-regular quaternionic functions $f:\cO\to\HH$ of class $C^1$. 
We find differential operators $T$ and $N$, with complex coefficients, such that a function $f$  is regular on $\OO$ if and only if $(N-jT)f=0$ on $\bO$ ($j$ a basic quaternion) and $f$ is harmonic on $\OO$.
As a consequence, by means of the identification of $\HH$ with $\CC^2$, we obtain a non-tangential holomorphicity condition which generalizes a result of Aronov and Kytmanov. We also show how the differential criterion and regularity are related to the $\db$-Neumann problem in $\CC^2$.

\end{abstract}

\maketitle
\thispagestyle{empty}

\section{Introduction}
\label{sec:Introduction}
Let $\OO$ be a smooth bounded domain in $\CC^2$. Let $\HH$ be the space of real quaternions $q=x_0+ix_1+jx_2+kx_3$, where  ${i,j,k}$ denote the basic quaternions. We identify $\HH$  with $\CC^2$ by means of the mapping that associates  the quaternion $q=z_1+z_2j$ with the pair $(z_1,z_2)=(x_0+ix_1,x_2+ix_3)$. 

In this paper we give a boundary differential criterion that characterizes \emph{(left) re\-gular} functions $f:\OO\to\HH$  (in the sense of Fueter) among harmonic functions. 

We show (Corollary \ref{cor1}) that there exist first order differential operators $T$ and $N$, with complex coefficients, such that a harmonic function $f:\OO\to\HH$, of class $C^1$ on $\cO$, is regular  if and only if $(N-jT)f=0$ on $\bO$. 

In order to obtain this result we study a related space of functions that satisfy a variant of the Cauchy-Riemann-Fueter equations, the space $\mathcal R(\OO)$ of \emph{$\psi$-regular} functions on $\OO$ (see \S\ref{sec:NotationsAndPreliminaries} for the precise definitions) for the particular choice $\psi=\{1,i,j,-k\}$ of the structural vector.   These functions have been studied by many authors (see for instance \cite{SV,MS,KS,No}). 
The space $\mathcal R(\OO)$ contains the identity mapping and any holomorphic mapping $(f_1,f_2)$ on $\OO$ defines a $\psi$-regular function $f=f_1+f_2j$. This is no more true if we replace the class of $\psi$-regular functions with that of regular functions. The definition of $\psi$-regularity is also equivalent to that of $q$-holomorphicity given by Joyce in \cite{J}, in the setting of hypercomplex manifolds. 

The space $\mathcal R(\OO)$ exhibits other interesting links with the theory of two complex variables. In particular, Vasilevski and  Shapiro \cite{VS} have shown that the Bochner-Martinelli kernel $U(\z,z)$ can be considered as a first complex component of the Cauchy-Fueter kernel associated to $\psi$-regular functions. This property was already observed by Fueter (see \cite F) in the general $n$-dimensional case, by means of an imbedding of $\CC^n$ in a real Clifford algebra.
Note that regular functions are in a simple correspondence with $\psi$-regular functions, since they can be obtained from them by means of a real coordinate reflection in $\HH$.

We prove (Theorem \ref{teo1}) that a harmonic function $f$ on $\OO$, of class $C^1$ on $\cO$, is $\psi$-regular  on $\OO$ if and only if $(\db_n-jL)f=0$ on $\bO$, where $\db_n$ is the normal part of $\db$ and $L$ is a tangential Cauchy-Riemann operator. 

This equation, which appeared in \cite{P} in connection with the characterization of the traces of pluriharmonic functions, can be considered as a generalization of both the CR-tangential equation $L(f)=0$  and the condition $\db_nf=0$ on $\bO$ that distinguishes holomorphic functions among complex harmonic functions  (Aronov and Kytmanov, see \cite{A,AK,KA}).

As an application of the differential condition for $\psi$-regular functions, we also obtain (Theorem \ref{teo4}) a differential criterion for holomorphicity of functions that generalizes, for a domain with connected boundary in $\CC^2$, the result of Aronov and Kytmanov.

In \S\ref{sec:DifferentialCriteriaForRegularityAndHolomorphicity} we give a weak formulation of the differential criterion of $\psi$-regularity, which makes sense also when the harmonic function $f$ is only continuous on the closure $\cO$. We obtain  trace theorems for $\psi$-regular functions, with applications to  holomorphic functions. Other results in the same vein were given by Pertici in \cite{Pe}, with generalizations to several quaternionic variables. 
Integral criteria for regularity were given also in \cite{SV,MS,VS}, under the assumption that the trace function satisfies a H\"older condition or belongs to a $L^p(\bO)$ space.

 In \S\ref{sec:RegularityAndTheDbNeumannProblem} we study the relation between regularity and the $\db$-Neumann problem in $\CC^2$ in the formulation given by Kytmanov in \cite K\S14--18. 
We are interested in a quaternionic analogue of the Hilbert transform, which relates one of the complex components of a \pr\ function to the boundary values of the other. 
We refer to \cite{KS} and \cite{RST} for generalizations of the Hilbert transform to the quaternionic setting. In these papers the functions considered are defined on plane or spatial domains, while we are interested in domains of two complex variables. In the latter case, pseudoconvexity becomes relevant, since such a domain is pseudoconvex if and only if every complex harmonic function on it is a complex component of a $\psi$-regular function (cf. Naser\cite{N} and N\=ono\cite{No1}). 

 In particular we show (Corollary \ref{cor4}) that if $\OO$ is a strongly pseudoconvex domain of class $C^\infty$ or a weakly pseudoconvex domain with  real-analytic boundary,  then the operator  that associates to $f=f_1+f_2j$ the restriction to $\bO$ of its first complex component $f_1$ induces an isomorphism between the quotient  spaces
${\mathcal R^\infty(\Omega)}/{A^\infty(\OO,\CC^2)}$ and ${C^\infty(\bO)}/{CR(\bO)}$, where $\mathcal R^\infty(\Omega)$ denotes the space of $\psi$-regular functions that are smooth up to the boundary and $A^\infty(\OO,\CC^2)$ is the space $Hol(\OO,\CC^2)\cap C^\infty(\cO,\CC^2)$

Some of the results contained in the present paper have been announced in \cite{P1}.
\section{Notations and preliminaries}
\label{sec:NotationsAndPreliminaries}

\subsection{}
\label{sec:NotationsAndPreliminaries1}
Let $\OO=\{z\in\Cn:\rho(z)<0\}$ be a bounded domain in $\Cn$ with boundary of
class $C^m, m\ge 1$. We assume $\rho\in C^m$ on $\CC^n$ and $d\rho\ne 0$ on
$\bO$. 

Let $\nu$ denote the outer unit normal to $\bO$ and $\tau=i\nu$. For every $f\in C^1(\cO)$,  we 
set $\db_nf=\frac12\left(\dd f\nu+i\dd f\tau\right)$ (see \cite {K}\S\S 3.3 and 14.2).

Then in a neighbourhood of
$\bO$ we have the decomposition of $\db f$ in the  tangential and the normal
parts \[\db f=\db_bf+\db_nf\frac{\db\rho}{\vert\db\rho\vert}.\] 

 The normal part
of $\db f$ on $\bO$ can also be expressed in the form \[\db_nf=\sum_k\dzk
f\dd\rho{z_k}\frac1{\vert\db\rho\vert},\]where $\vert\db\rho\vert^2= \sum_{k=1}^n\left\vert\dd\rho{\bar z_k}\right\vert^2$, or,
by means of the Hodge $*$-operator
and the Lebesgue surface measure $d\sigma$, as $\db_nf d\sigma={*\db f
}_{\vert_{\bO}}$. 

\subsection{}
\label{sec:NotationsAndPreliminaries2}
We recall the definition of  tangential Cauchy-Riemann operators (see for example \cite{R}\S 18). A linear first-order differential operator $L$ is \emph{ tangential} to $\bO$ if $(L\rho)(\zeta)=0$ for each point $\zeta\in\bO$. A tangential operator of the form \[L=\sum_{j=1}^n a_j\dd{}{\bar z_j}\] is called a \emph{ tangential Cauchy-Riemann} operator.
The operators
\[L_{jk}=\frac1{\vert \db\rho\vert}\left (\dd\rho{\bar z_k}\dd{}{\bar z_j}-\dd\rho{\bar z_j}\dd{}{\bar z_k}\right )
\textrm{ ,\quad }1\le j<k\le n,\] are tangential and the corresponding vectors at $\zeta\in\bO$ span (not independently when $n>2$) the (conjugate) complex tangent space to $\bO$ at $\zeta$. Then a function $f\in C^1(\bO)$ is a CR function if and only if $L_{jk}(f)=0$ on $\bO$ for every $j,k$, or, equivalently, if $L(f)=0$ for each tangential Cauchy-Riemann operator $L$. In particular, when $n=2$, $f$ is a CR function if and only if $L(f)=0$ on $\bO$.

\subsection{}
\label{sec:NotationsAndPreliminaries3}
Let $\HH$ be the algebra of quaternions.  The elements of $\HH$ have the form \[q=x_0+ix_1+jx_2+kx_3,\] where $x_0,x_1,x_2,x_3$ are real numbers and ${i,j,k}$ denote the basic quaternions. \par
We identify   the space $\CC^2$ with the set $\HH$ by means of the mapping that associates  the pair $(z_1,z_2)=(x_0+ix_1,x_2+ix_3)$ with the quaternion $q=z_1+z_2j$.
The commutation rule is then $aj=j\bar a$ for every $a\in\CC$, and the quaternionic conjugation is \[\bar q=x_0-ix_1-jx_2-kx_3=\bar z_1-z_2j.\]
We refer to \cite{S} and  \cite{BDS} for the theory of quaternionic analysis and its generalization represented by Clifford analysis. 
We will denote by $\mathcal D$ the left Cauchy-Riemann-Fueter operator \[\mathcal D=\dd{}{x_0}+i\dd{}{x_1}+j\dd{}{x_2}+k\dd{}{x_3}.\]

A quaternionic $C^1$ function $f=f_1+f_2j$, is \emph{ (left-)regular} on a domain $\OO\subseteq\HH$ if $\mathcal Df=0$ on $\OO$.
We prefer to work with another class of regular functions, defined by the Cauchy-Riemann-Fueter operator associated with the structural vector $\psi=\{1,i,j,-k\}$:
\[\mathcal D'=\dd{}{x_0}+i\dd{}{x_1}+j\dd{}{x_2}-k\dd{}{x_3}=2\left(\dd{}{\bar z_1}+j\dd{}{\bar z_2}\right).\]
A quaternionic $C^1$ function $f=f_1+f_2j$, is called \emph{ (left-)$\psi$-regular} on a domain $\OO$, if $\mathcal D'f=0$ on $\OO$. This condition is equivalent to the following system of complex differential equations:
\begin{equation}\label{eq1}
 \dd{f_1}{\bar z_1}=\dd{\overline {f_2}}{z_2},\quad \dd {f_1}{\bar z_2}=-\dd{\overline {f_2}}{z_1}
\end{equation}
or to the equation $*\db f_1=-\frac12\partial(\overline{f_2} d\bar z_1\wedge d\bar z_2)$.  
Note that the identity mapping is \pr, and any holomorphic mapping $(f_1,f_2)$ on $\OO$ defines a $\psi$-regular function $f=f_1+f_2j$. This is no more true if we replace the class of $\psi$-regular functions with that of regular functions.

If $\OO$ is pseudoconvex, every complex harmonic function $f_1$ on $\OO$ is a complex component of a $\psi$-regular function $f=f_1+f_2j$, since the $(1,2)$-form $*\db f_1$ is $\partial$-closed on $\OO$. (cf. \cite{N} and \cite{No1} for  this result and its converse). Regular and $\psi$-regular functions are real analytic on $\OO$, and they are harmonic with respect to the Laplace operator in $\RR^4$. 
We refer, for instance, to \cite{SV,MS,KS,No} for the properties of structural vectors and $\psi$-regular functions (in these references they are called $\psi$-hyperholomorphic functions). 

\begin{remark}Let $\gamma$ be the transformation of $\CC^2$ defined by $\gamma(z_1,z_2)=(z_1,\bar z_2)$. Then a $C^1$ function $f$ is regular on the domain $\OO$ if, and only if, $f\circ \gamma$ is $\psi$-regular on $\gamma^{-1}(\OO)$.
\end{remark}

\subsection{}
\label{sec:NotationsAndPreliminaries4}
A definition equivalent to $\psi$-regularity has been given by Joyce in \cite{J}, in the setting of hypercomplex manifolds. Joyce introduced the module of \emph{q-holomorphic} functions  on a hypercomplex manifold. A hypercomplex structure on the manifold  $\HH$ is given by the complex structures $J_1,J_2$ defined on $T^*\HH\simeq\HH$  by left multiplication by $i$ and $j$, and by $J_3=J_1J_2$. These structures act on the basis elements in the following way:
\begin{eqnarray*}J_1dx_0=-dx_1,\  J_1dx_2&=&-dx_3,\  J_2dx_0=-dx_2,\ J_2dx_1=dx_3,\\
  J_3dx_0&=&dx_3,\ J_3dx_1=dx_2.
\end{eqnarray*}

An easy computation shows that a differentiable function $f:\OO\rightarrow \HH$ is $\psi$-regular on $\OO$ if and only if it is q-holomorphic, that is
\[df+iJ_1(df)+jJ_2(df)+kJ_3(df)=0\quad \textrm{ on }\OO.\]
Equivalently, if $f=f^0+if^1+jf^2+kf^3$ is the real decomposition of $f$ with respect to the standard basis, $f$ is $\psi$-regular if and only if the real equation
\[df^0=J_1(df^1)+J_2(df^2)+J_3(df^3)\]
is satisfied on $\OO$. Returning to complex components, we can rewrite equations (\ref{eq1}) by means of the complex structure $J_2$ as follows:
\[\db f_1=J_2(\partial\overline f_2).\]

Note that our notations are slightly different from those of Joyce \cite{J}. With our choices, $\psi$-regular functions form a \emph{right} $\HH$-module and a complex valued function $f=f^0+if^1:\OO\rightarrow\CC$ is holomorphic w.r.t. a complex structure $J$ if it satisfies the Cauchy-Riemann equations $df^0=J(df^1)$ on $\OO$ or, equivalently, $df+iJ(df)=0$. 

Let $J_p$ be the complex structure defined by an imaginary quaternion $p=p_1i+p_2j+p_3k$ in the unit sphere $S^2$ as $J_p=p_1J_1+p_2J_2+p_3J_3$. Then every $J_p$-holomorphic function $f=f^0+if^1:\OO\rightarrow\CC$ defines a $\psi$-regular function $\tilde f=f^0+pf^1$ on $\OO$. The original function $f$ can be recovered from $\tilde f$ by means of the formula
\[f=\textrm{ Re}(\tilde f)+i\textrm{ Re}(-p\tilde f).\]

The \pr\ functions obtained in this way can be identified with the holomorphic functions from the complex manifold $(\OO,J_p)$ to the manifold $(\CC_p,L_p)$, where $\CC_p=\langle 1,p\rangle$ is a copy of $\CC$ in $\HH$ and $L_p$ is the complex structure defined on $T^*\CC_p\simeq\CC_p$  by left multiplication by $p$. In fact, 
\[d\tilde f=J_p(df^1)-pJ_p(df^0)=-pJ_p(df^0+pdf^1)=-pJ_p(d\tilde f).\]

More generally, holomorphic maps w.r.t.\ any complex structure $J_p$ induce $\psi$-regular functions, since for any positive orthonormal basis $\{1,p,q,pq\}$ of $\HH$ ($p,q\in S^2$), $\psi$-regular functions are the solutions of the equation
\[\db_p f_1=J_q(\partial_p\overline f_2),\] where $f=(f^0+pf^1)+(f^2+pf^3)q=f_1+f_2q$, $\overline f_2=f^2-pf^3$  and $\db_p$ is the Cauchy-Riemann operator w.r.t.\ $J_p$:
\[\db_p=\frac12\left(d+pJ_p\circ d\right).\] 
We shall return in a subsequent paper on the problem of the characterization of the (proper) submodule of holomorphic functions in the $\HH$-module of $\psi$-regular functions.

\subsection{}
\label{sec:NotationsAndPreliminaries5}
Let's denote by $G$ the Cauchy-Fueter quaternionic kernel defined by \[G(p-q)=\frac1{2\pi^2}\frac{\bar p-\bar q}{|p-q|^4},\] and by $G'$ the Cauchy kernel for $\psi$-regular functions:
\[G'(p-q)=\frac1{2\pi^2}\frac{y_0-x_0-i(y_1-x_1)-j(y_2-x_2)+k(y_3-x_3)}{|p-q|^4},\]
where $p=y_0+iy_1+jy_2+ky_3$, $q=x_0+ix_1+jx_2+kx_3$.

Let $\sigma(q)$ be the  quaternionic valued 3-form \[\sigma(q)=dx[0]-idx[1]+jdx[2]+kdx[3],\] where $dx[k]$ denotes the product of $dx_0,dx_1,dx_2,dx_3$ with $dx_k$ deleted. Then the Cauchy-Fueter integral formula for left-$\psi$-regular functions on $\OO$ that are continuous on $\cO$, holds true:
\[\int_{\bO}G'(p-q)\sigma(p)f(p)=
\left\{
\begin{array}{c}
	f(q)\textrm{\quad for\ }q\in\OO,\\0\textrm{\quad for\ }q\notin\cO.
\end{array}
\right.
\]

In \cite{VS}  (see also \cite F and \cite{MS}) it was shown that, for a family of structural vectors, including $\{1,i,j,-k\}$, the two-dimensional \BM\ 
form $U(\z,z)$ can be considered as a first complex component of the Cauchy-Fueter kernel associated to $\psi$-regular functions. Let $q=z_1+z_2j$, $p=\z_1+\z_2j$. Then
\[G'(p-q)\sigma(p)=U(\z,z)+\omega(\z,z) j,\] where $\omega(\z,z)$ is the following complex $(1,2)$-form:
\[\omega(\z,z)=-\frac1{4\pi^2|\z-z|^4}\left((\bar\z_1-\bar z_1)d\z_1+(\bar\z_2-\bar z_2)d\z_2\right)\wedge \overline{d\z}.\]
Here ${d\z}={d\z_1}\wedge{d\z_2}$ and we choose the orientation of $\CC^2$ given by the volume form $\frac14dz_1\wedge dz_2\wedge\overline{dz_1}\wedge\overline{dz_2}$.

\section{Differential criteria for regularity and holomorphicity}
\label{sec:DifferentialCriteriaForRegularityAndHolomorphicity}

\subsection{}
We now rewrite the representation formula of Cauchy-Fueter for $\psi$-regular functions in complex form. We use results from \cite{MS} to relate the form $\omega(\z,z)$ to the tangential operator $L_{12}$, that we will denote simply by $L$.  We show that the \BM\ formula can then be applied to obtain a criterion that distinguishes regular functions among harmonic functions on a domain $\OO$ in $\CC^2=\HH$.

Let $g(\z,z)=\frac1{4\pi^2}|\z-z|^{-2}$ be the fundamental solution of the complex laplacian on $\CC^2$. 

\begin{proposition}
\label{pro1}
Let $\OO$ be a bounded domain of class $C^1$ in $\HH$. Let $f:\OO\to\HH$ be a quater\-nio\-nic function, of class $C^1$ on $\cO$. Then $f$ is (left-)$\psi$-regular  on $\OO$ if, and only if, the following representation formula holds on $\OO$:
\[f(z)=\int_{\bO}U(\z,z)f(\z)+2\int_{\bO}g(\z,z)jL(f(\z))d\sigma
\]where $d\sigma$ is the Lebesgue measure on $\bO$ and the tangential operator $L$ acts on $f=f_1+f_2j$ as $L(f)=L(f_1)+L(f_2)j$.
\end{proposition}
\begin{proof}The integral  of Cauchy-Fueter in complex form is \[\int_{\bO}G'(p-q)\sigma(p)f(p)=\int_{\bO}U(\z,z)f(\z)+\int_{\bO}\omega(\z,z) jf(\z).\]
From Proposition 6.3 in \cite{MS}, we get that the last integral is equal to 
\begin{eqnarray*}
\int_{\bO}\omega(\z,z)\overline{f_1}j-\int_{\bO}\omega(\z,z)\overline{f_2}&=&
2\int_{\bO}g(\z,z)\left(\overline{L(f_1)}j-\overline{L(f_2)}\right)d\sigma\\&=&2\int_{\bO}g(\z,z)jL(f)d\sigma.
\end{eqnarray*}
Then the result follows from the Cauchy-Fueter integral formula for $\psi$-regular functions.
\end{proof}

If $f=f_1+f_2j$ is a $\psi$-regular function on $\OO$, of class $C^1$ on $\cO$, then from equations (\ref{eq1}) we get that it satisfies the equation \begin{equation}\label{eq2}
(\db_n-jL)f=0\textrm{\quad on }\bO,
\end{equation}
since $\db_nf_1=-\overline {L(f_2)}$, $\db_nf_2=\overline{L(f_1)}$ on $\bO$.\par
This equation was introduced in \cite{P}\S4 in connection with the characterization of the traces of pluriharmonic functions. It can be considered as a generalization both of the CR-tangential equation $L(f)=0$ (for a complex-valued $f$) and of the condition $\db_nf=0$ on $\bO$ that distinguishes holomorphic functions among complex harmonic functions on $\OO$ (what is called the homogeneous $\db$-Neumann problem for functions, see \cite{AK} and \cite{K}\S15).

\begin{theorem}\label{teo1}
Let $\OO$ be a bounded domain  in $\HH$, with boundary of class $C^1$. Let $f=f_1+f_2j:\OO\to\HH$ be a harmonic function on $\OO$, of class $C^1$ on $\cO$. Then, $f$ is (left-)$\psi$-regular  on $\OO$ if, and only if, \[(\db_n-jL)f=0\textrm{\quad on }\bO.\]
\end{theorem}
\begin{proof}
It remains to prove the sufficiency of condition (\ref{eq2}) for $\psi$-regularity of harmonic functions. For every $z\in\OO$,  it follows from  the \BM\ integral representation for complex harmonic functions on $\OO$ (see for example \cite{K}\S1.1), that $f(z)=f_1(z)+f_2(z)j$ is equal to
\begin{eqnarray*}
\int_{\bO}U(\z,z)f_1(\z)+2\int_{\bO}g(\z,z)\db_n f_1(\z)d\sigma&&\\+
&&\left(\int_{\bO}U(\z,z)f_2(\z)+2\int_{\bO}g(\z,z)\db_n f_2(\z)d\sigma\right)j.
\end{eqnarray*} 
If $\db_n f=jL(f)$ on $\bO$, then we obtain 
\begin{eqnarray*}
f(z)&=&\int_{\bO}U(\z,z)f(\z)+2\int_{\bO}g(\z,z)\db_n f(\z)d\sigma\\&=&\int_{\bO}U(\z,z)f(\z)+2\int_{\bO}g(\z,z)jL(f(\z))d\sigma.
\end{eqnarray*}
The result now follows from Proposition \ref{pro1}.
\end{proof}

Let $N$ and $T$ be the differential operators, defined in a neighbourhood of $\bO$, as
\[N=\dd\rho{z_1}\dd{}{\bar z_1}+\dd\rho{\bar z_2}\dd{}{z_2},\quad
T=\dd\rho{z_2}\dd{}{\bar z_1}-\dd\rho{\bar z_1}\dd{}{z_2}.\]
$T$ is a tangential (not Cauchy-Riemann) operator to $\bO$, while $N$ is non-tangential, such that $N(\rho)={\vert \db\rho\vert}^2$, $\textrm{Re}(N)={\vert \db\rho\vert}\textrm{Re}(\db_n)$. The remark made at the end of \S\ref{sec:NotationsAndPreliminaries3} shows that Theorem \ref{teo1} gives also a boundary condition for regularity of a harmonic function on $\OO$.

\begin{corollary}\label{cor1}
Let $\OO$ be a $C^1$-bounded domain  in $\HH$. Let $f=f_1+f_2j:\OO\to\HH$ be a harmonic function on $\OO$, of class $C^1$ on $\cO$. Then, $f$ is (left-)regular  on $\OO$ if, and only if, \[(N-jT)f=0\textrm{\quad on }\bO.\]
\end{corollary}

\subsection{}
We now give a weak formulation of the differential criterion of $\psi$-regularity, which makes sense for example when the harmonic function $f$ is only continuous on the closure $\cO$. \par Let $Harm^1(\cO)$ denote the space of complex harmonic functions on $\OO$, of class $C^1$ on $\cO$.
By application of the Stokes' Theorem, of the complex Green formula \[\int_{\bO}g\db_n hd\sigma=\int_{\bO}h\partial_n g d\sigma\textrm{\quad }\forall\ g,h\in Harm^1(\cO),\] and of the equality $\db f\wedge d\z_{|\bO}=2L(f)d\sigma$ on $\bO$,
we see that the equations $\db_nf_1=-\overline {L(f_2)}$, $\db_nf_2=\overline{L(f_1)}$ on $\bO$, imply the following (complex) integral conditions: for every  function $\phi\in Harm^1(\cO)$,
\begin{equation}\label{eq3}
\int_{\bO}\overline{f_1}{*\db\phi}=\frac12\int_{\bO}f_2\db(\phi d\z),\quad \int_{\bO}\overline{f_2}{*\db\phi}=-\frac12\int_{\bO}f_1\db(\phi d\z).
\end{equation}
These are equivalent  to one quaternionic condition, which is then necessary for the $\psi$-regularity of $f\in C^1(\cO)$:
\[\int_{\bO}\bar f \left({*\db\phi} -\frac12 j\db(\phi d\z)\right)=0\quad \forall\ \phi\in Harm^1(\cO),\]
that can be rewritten also as: 
\begin{equation}\label{eq4}
\int_{\bO}\bar f \left(\db_n-jL\right)(\phi)\ d\sigma=0\quad \forall\ \phi\in Harm^1(\cO).
\end{equation}

Now we consider the sufficiency of the integral condition (\ref{eq4}) when $f$ is only continuous on $\cO$.

\begin{theorem}\label{teo2}
Let $\OO$ be a bounded domain  in $\HH$, with boundary of class $C^1$. Let $f:\bO\to\HH$ be a continuous function. Then, there exists a (left-)$\psi$-regular function $F$ on $\OO$, continuous on $\cO$, such that $F_{|\bO}=f$, if and only if $f$ satisfies the condition (\ref{eq4}).
\end{theorem}
\begin{proof}Let $F^+$ and $F^-$ be the $\psi$-regular functions defined respectively on $\OO$ and on $\CC^2\setminus\cO$ by the Cauchy-Fueter integral of $f$: \[F^\pm (z)=\int_{\bO}U(\z,z)f(\z)+\int_{\bO}\omega(\z,z) jf(\z).\]From the equalities $U(\z,z)=-2{*\partial_\z g(\z,z)}$, $\omega(\z,z)=-\partial_\z (g(\z,z) d\bar\z)$, we get that \[\overline{F^-(z)}=-2\int_{\bO}\overline{f(\z)}{*\db_\z{g(\z,z)}}+\int_{\bO}\overline{f(\z)} j\db_\z({g(\z,z)} d\z)\]for every $z\notin\cO$.
If (\ref{eq4}) holds, then $F^-$ vanishes identically on $\CC^2\setminus\OO$. 
As in the complex variable case, extended to the quaternionic case in Lemma 3 of \cite{Pe}, this implies 
that also $F^+$ has a continuous extension on $\cO$, and $F^+=f$ on $\bO$. Conversely, if $F\in C(\cO)$ is a $\psi$-regular function on $\OO$ with trace $f$ on $\bO$, and $\OO_\epsilon=\{z\in\OO:\rho<\epsilon\}$, then $F$  satisfies (\ref{eq4}) on $\bO_\epsilon$ for every small $\epsilon<0$. Passing to the limit as $\epsilon\to 0$, we obtain (\ref{eq4}).
\end{proof}
\begin{remark}If $f$ satisfies a H\"older condition on $\bO$, the same result can be obtained using the Sokhotski-Plemelj formula (see \cite{SV}\S3.6).
In \cite{SV,MS,VS} other similar integral criteria were given, assuming that the trace function belongs to a H\"older or a $L^p(\bO)$ class.
\end{remark}
\begin{remark}
In the orthogonality condition (\ref{eq4}) it is sufficient to consider functions $\phi\in Harm^1(\cO)$ that are of class $C^\infty$ on a neighbourhood of $\cO$.
\end{remark}

\subsection{}
In Theorem \ref{teo2},  the boundary of $\OO$ is not required to be connected. If $\bO$ is connected, we can improve the result and show that only one of the complex conditions (\ref{eq3}) is sufficient for the $\psi$-regularity of  the harmonic extension of $f$.

\begin{theorem}\label{teo3}
Let $\OO$ be a bounded domain  in $\HH$, with connected boundary $\bO$ of class $C^1$. Let $f:\bO\to\HH$ be a continuous function. Then, if $f$ satisfies one of the conditions (\ref{eq3}), there exists a (left-)$\psi$-regular function $F$ on $\OO$, continuous on $\cO$, such that $F_{|\bO}=f$.
\end{theorem}
\begin{proof}Assume that \[ \int_{\bO}\overline{f_2}{*\db\phi}=-\frac12\int_{\bO}f_1\db(\phi d\z)\quad \forall\ \phi\in Harm^1(\cO).\]
We use the same notation as in the proof of Theorem \ref{teo2}. We get that 
\begin{eqnarray*}
\overline{F^-(z)}=-2\int_{\bO}\overline{f(\z)}{*\db_\z{g(\z,z)}}+\int_{\bO}\overline{f(\z)} j\db_\z({g(\z,z)} d\z)\\=-2\int_{\bO}\overline{f_1(\z)}{*\db_\z{g(\z,z)}}+\int_{\bO}f_2(\z)\db_\z({g(\z,z)})
\end{eqnarray*}
for every $z\notin\cO$. Therefore, $F^-$ is a complex-valued, $\psi$-regular function on $\CC^2\setminus\cO$. The system of equations (\ref{eq1}) then implies that $F^-$ is a holomorphic function. Since $\bO$ is connected, from Hartogs' Theorem it follows that $F^-$ can be holomorphically continued to the whole space. Let $\tilde F^-$ be such extension. Then $F=F^+-{\tilde F^-}_{|\OO}$ is a $\psi$-regular function on $\OO$, continuous on $\cO$, such that $F_{|\bO}=f$. If the first condition in (\ref{eq3}) is satisfied, it is sufficient to consider the function $fj=-f_2+f_1j$ in place of $f$.
\end{proof}

\begin{remark}
The hypothesis on $f$ in the preceding theorem is satisfied, for example, when $f$ is of class $C^1$ on $\cO$, harmonic on $\OO$, and one of the equations $\db_nf_1=-\overline {L(f_2)}$, $\db_nf_2=\overline{L(f_1)}$ holds on $\bO$.
\end{remark}
\begin{remark}
The connectedness of $\bO$ is a necessary condition in Theorem \ref{teo3}: consider a locally constant function on $\bO$. For example, il $f_2=0$ on $\OO$ and $f_1$ takes two distinct values on the components of $\bO$, then $\db_nf_2=\overline{L(f_1)}$ on $\bO$, but $\db_nf_1\ne 0$ on $\bO$, since otherwise $f_1$ would be holomorphic on $\OO$.
\end{remark}

The preceding result can be easily generalized in the following form: 
\begin{corollary}
\label{cor2}
Let $\OO$ be as above. Let $a,b\in\CC$ be two complex numbers such that $(a,b)\ne (0,0)$.\par
(i) If $f$ is of class $C^1$ on $\cO$, harmonic on $\OO$, and such that \[a\ \db_n f_1+b\ \db_n f_2=-a\ \overline {L(f_2)} +b\ \overline {L(f_1)}\textrm{\quad on\ }\bO,\]then $f$ is $\psi$-regular on $\OO$.\par
(ii) If $f$ is continuous on $\bO$, such that \[\int_{\bO}(a\overline{f_1}+b\overline{f_2}){*\db\phi}=\frac12\int_{\bO}(af_2-bf_1)\db(\phi d\z)\quad \forall\ \phi\in Harm^1(\cO),\] then there exists a $\psi$-regular function $F$ on $\OO$, continuous on $\cO$, such that $F_{|\bO}=f$.
\end{corollary}

Theorem \ref{teo3} can be applied,   in the case of connected boundary  in $\CC^2$,  to obtain the following result of Aronov and Kytmanov (cf. \cite {A,AK,KA}), which holds in $\CC^n$, $n>1$: if $f$ is a complex harmonic function on $\OO$, of class $C^1$ on $\cO$, such that $\db_n f=0$ on $\bO$, then $f$ is holomorphic. It is sufficient to take $f_1=f,f_2=0$.

More generally, we can deduce a differential criterion for holomorphicity of functions on a domain with connected boundary in $\CC^2$, analogous to those proposed in \cite{C} and investigated  in \cite{K1} and \cite{K}\S23.2.

\begin{theorem}
\label{teo4}Let $\OO$ be a bounded domain  in $\CC^2$, with connected boundary of class $C^1$. Let $h=(h_1,h_2):\OO\to\CC^2$ be a holomorphic mapping of class $C^1$ on $\cO$, such that $h(\z)\ne0$ for every $\z\in\bO$.\par (i) If $h_1f,h_2f \in Harm^1(\cO)$ and $f:\cO\to\CC$ satisfies the differential condition \[h_1\db_nf=\overline{h_2L(f)}\textrm{\quad on }\bO,\]then $f$ is holomorphic on $\OO$.\par
(ii) If $h_1f$ and $h_2f$ are harmonic on $\OO$, $f:\cO\to\CC$ is continuous and it satisfies the integral condition
\[\int_{\bO}\overline{h_1f}{*\db\phi}=\int_{\bO}h_2f\db(\phi d\z)\quad \forall\ \phi\in Harm^1(\cO),\]then $f$ is holomorphic on $\OO$.
\end{theorem}
\begin{proof}
It is sufficient to prove (ii). Let $h'_2=-2h_2$. From Theorem \ref{teo3} we get that the harmonic function $h'_2f+h_1fj$ is $\psi$-regular on $\OO$. Let $\OO_\epsilon=\{z\in\OO:\rho<\epsilon\}$ for a small $\epsilon <0$. Then  the following equalities hold on $\bO_\epsilon$: \[\db_n(h_1f)=\overline{L(h'_2f)},\quad \db_n(h'_2f)=-\overline{L(h_1f)}.\]From the holomorphicity of $h$, we then obtain \[h_1\db_n(f)=\bar h'_2\overline{L(f)},\quad h'_2\db_n(f)=-\bar h_1\overline{L(f)}\quad \textrm {on }\bO_\epsilon,\]which implies $\db_nf=L(f)=0$ on $\bO_\epsilon$ for every $\epsilon$ sufficiently small, such that $h\ne0$ on $\bO_\epsilon$. This means that there exists a holomorphic extension $F_\epsilon$ of $f$  on $\OO_\epsilon$. From the equality of the harmonic functions $h_jF_\epsilon=h_jf$ on $\OO_\epsilon$, for $j=1,2$, we get $F_\epsilon=f_{|\OO_\epsilon}$. Then $f$ is holomorphic on the whole domain $\OO$.
\end{proof}

\section{Regularity and the $\db$-Neumann problem}
\label{sec:RegularityAndTheDbNeumannProblem}

\subsection{}
We recalled before the following fact proved by Naser and N\=ono \cite{N,No1}: $\OO$ is pseudoconvex if and only if every complex harmonic function on $\OO$ is a complex component of a $\psi$-regular function.
Now we are interested in the boundary values of $\psi$-regular functions and in the quaternionic analogue of the Hilbert transform. We want to express one of the complex components of a \pr\ function by means of the other (the two components are then a pair of \lq conjugate harmonic\rq\ functions). We show how this is related to the $\db$-Neumann problem in $\CC^2$.

We refer, for instance, to \cite{KS} and \cite{RST} for generalizations of the Hilbert transform to the quaternionic setting. In these references the \pr\ (or more generally, $(\psi,\alpha)$-hyperholomorphic) functions are defined on plane or spatial domains. Here we are interested in domains of two complex variables where, as shown before, pseudoconvexity becomes relevant. We refer to \cite{BDS1} for other multidimensional generalizations of the concept of conjugate harmonic functions from the complex case to the Clifford space case.

The $\db$-Neumann problem for complex functions 
$\square f=(\db^*\db+\db\db^*)f=\psi$ in $\OO$, $\db_n f=0$ on $\bO$ is
equivalent,  in the smooth case, to the problem \[\db_n g=\phi\textrm{ on }
\bO\textrm{,\quad}g\textrm{ harmonic in }\OO\] (see \cite {K}\S 14). The compatibility condition for this problem is
\begin{equation}\label{eq5}
\int_{\bO}\phi\overline hd\sigma=0
\end{equation}
for every $h$ holomorphic in a
neighbourhood of $\cO$. We now use the solvability of this problem in strongly pseudoconvex domains of $\CC^2$ to obtain some results on regular functions.

Let $\OO$ be a bounded domain in $\CC^2$ with connected, $C^\infty$-smooth boundary. We denote by $W^s(\OO)$ ($s\ge1$) the complex Sobolev space,  and by $G^s(\OO)$ the space of harmonic functions in $W^s(\OO)$. The space $G^s(\OO)$ is isomorphic to  $W^{s-1/2}(\bO)$ through the restriction operator.

\begin{theorem}
\label{teo5}
Let $\OO$ be a bounded strongly pseudoconvex domain in $\CC^2$ with connected boundary of class $C^\infty$. Let $f_1\in W^{s-1/2}(\bO)$, where $s\ge3$. We identify $f_1$ with its harmonic extension in $G^s(\OO)$. Then there exists a function $f_2\in G^{s-2}(\OO)$ such that $f=f_1+f_2j$ is a $\psi$-regular function on $\OO$. 
\end{theorem}
\begin{proof}We show that the function $\phi=\overline{L(f_1)}\in G^{s-1}(\OO)$ satisfies the compatibility condition (\ref{eq5}). If $h$ is holomorphic in a
neighbourhood of $\cO$, \[\int_{\bO}L(f_1)hd\sigma=\frac12\int_{\bO}h\db(f_1dz)=0,\]since $h$ is a CR function on $\bO$. Then we can apply a result of Kytmanov \cite{K}\S18.2 and get a solution $f_2\in G^{s-2}(\OO)$ of the $\db$-Neumann problem $\db_nf_2=\overline{L(f_1)}$ on $\bO$. If $s\ge5$, then $f=f_1+f_2j$ is continuous on $\cO$ by Sobolev embedding. From Theorem \ref{teo3} we get that $f$ is $\psi$-regular on $\OO$, since it satisfies the second condition in (\ref{eq3}). In any case, $f_2\in L^2(\bO)$ since $s\ge3$. Then the result follows from the $L^2(\bO)$-version of Theorem \ref{teo3}, that can be proved as before using the results in 
\cite{SV}\S3.7. 
\end{proof}
\begin{remark}There is a unique solution $f_2$ of $\db_nf_2=\overline{L(f_1)}$ on $\bO$ that is orthogonal to holomorphic functions in $L^2(\bO)$. It is given by the bounded Neumann operator $N_{\OO}$: $f_2=N_\OO(\overline{L(f_1)})$.
\end{remark}
\begin{corollary}
\label{cor3}
Suppose $\OO$ is a bounded strongly pseudoconvex domain in $\CC^2$ with connected boundary of class $C^\infty$. Let $f_1:\bO\to\CC$ be of class $C^\infty$. Then there exists a $\psi$-regular function $f$ on $\OO$, of class $C^\infty$ on $\cO$, such that the first complex component of the restriction $f_{|\bO}$ to $\bO$ is $f_1$.
\end{corollary}
\begin{remark}
The preceding statements remain true if $\OO$ is a bounded weakly pseudoconvex domain in $\CC^2$ with connected real-analytic boundary, since on these domains the $\db$-Neumann  problem for smooth functions is solvable (cf. \cite{K}\S18).
\end{remark}

\subsection{}
We denote by $\mathcal R(\OO)$ the right $\HH$-module of (left-)$\psi$-regular functions on $\OO$ and by $\mathcal R^\infty(\Omega)$ the functions in $\mathcal R(\OO)$ that are of class $C^\infty$ on $\cO$. We consider the space of holomorphic maps $Hol(\OO,\CC^2)$ as a real subspace of $\mathcal R(\OO)$ by identification of the map $(f_1,f_2)$ with $f=f_1+f_2j$.

If $\OO$ is pseudoconvex, it follows from what observed in \S\ref{sec:NotationsAndPreliminaries3} that the map that associates to $f=f_1+f_2j$ the first complex component $f_1$ induces an isomorphism between the quotient real spaces ${\mathcal R(\OO)}/{Hol(\OO,\CC^2)}$ and  ${Harm(\OO)}/{\mathcal O(\OO)}.$ 

Now we are also interested in the regularity up to the boundary of the functions. Let $A^\infty(\OO,\CC^2)=Hol(\OO,\CC^2)\cap C^\infty(\cO,\CC^2)$ be identified with a $\RR$-subspace of 
$\mathcal R^\infty(\Omega)$.

Let $C:\mathcal R^\infty(\Omega)\to C^\infty(\bO)$ be the linear operator that associates to $f=f_1+f_2j$ the restriction to $\bO$ of its first complex component $f_1$. From the Corollary \ref{cor3} and the remark preceding it, we get a right inverse $R$ of $C$. The function $R(f_1)$ is uniquely determined by the orthogonality condition with respect to the functions holomorphic on a neighbourhood of $\cO$:
\[\int_{\bO}(R(f_1)-f_1)\overline h d\sigma=0\quad\forall h\in \mathcal O(\cO).\]
Note that $f_1\in CR(\bO)\cap C^\infty(\bO)$ if and only if $R(f_1)\in A^\infty(\OO,\CC^2)$. Besides, the operator $C$ has kernel \[\ker C=\{f_2j : f_2\in A^\infty(\OO)\}=A^\infty(\OO)j\] where $A^\infty(\OO)$ is the space of the holomorphic functions on $\OO$ that are $C^\infty$ up to the boundary. Then $C$ induces the following isomorphism of real spaces:
\[\tilde C\ :\ \frac{\mathcal R^\infty(\Omega)}{ A^{\infty}(\OO)j}\to C^\infty(\bO).\]

\begin{corollary}
\label{cor4}
Let $\OO$ be a bounded  domain in $\CC^2$ with connected boundary. Suppose that $\OO$ is a strongly pseudoconvex domain of class $C^\infty$ or a weakly pseudoconvex domain with  real-analytic boundary. Then the operator $C$ induces an isomorphism of real spaces:
\[\hat C\ :\ \frac{\mathcal R^\infty(\Omega)}{A^\infty(\OO,\CC^2)}\to \frac{C^\infty(\bO)}{CR(\bO)}.\]
\end{corollary}

\begin{remark}If $\OO$ is a $C^\infty$-smooth bounded pseudoconvex domain, from application of Kohn's Theorem on the solvability of the $\db$-problem  to the equation $*\db f_1=-\frac12\partial(\overline{f_2} d\bar z_1\wedge d\bar z_2)$ we can still deduce the isomorphism of Corollary \ref{cor4}.
\end{remark}

\end{document}